\documentclass[reqno,centertags,12pt]{amsart}
\usepackage{amsmath,amsthm,amscd,amssymb,latexsym,verbatim}

\usepackage{xcolor,graphicx,epsf,cite}

\usepackage{tikz}
\pgfdeclarelayer{nodelayer}
\pgfdeclarelayer{edgelayer}
\pgfsetlayers{edgelayer,nodelayer,main}
\tikzstyle{arrow}=[draw=black,arrows=-latex]

-\marginparwidth \oddsidemargin -4mm \evensidemargin \oddsidemargin

\newtheorem{theorem}{Theorem}[section]

\newtheorem{lemma}[theorem]{Lemma}

\theoremstyle{definition}

\theoremstyle{remark}

\newcounter{smalllist}

\DeclareMathOperator*{\sgn}{sgn}

\allowdisplaybreaks
\numberwithin{equation}{section}

\newcommand{\lb}{\label}

\newcommand{\beq}{\begin{equation}}
\newcommand{\eeq}{\end{equation}}

\newcommand{\bal}{\begin{align}}
\newcommand{\eal}{\end{align}}
\newcommand{\bals}{\begin{align*}}
\newcommand{\eals}{\end{align*}}

\newcommand{\bbR}{{\mathbb{R}}}

\newcommand{\bbT}{{\mathbb{T}}}

\newcommand{\eps}{\varepsilon}

\newcommand{\til}{\tilde}

\begin{document}
\title[Double-Exponential Growth for Euler on the Half-Plane]
{Maximal Double-Exponential Growth \\ for the Euler Equation on the Half-Plane}

\author{Andrej Zlato\v s}
\address{\noindent Department of Mathematics \\ UC San Diego \\ La Jolla, CA 92093, USA \newline Email:
zlatos@ucsd.edu}


\begin{abstract} 
We show that smooth  solutions to the Euler equation on the half-plane  can exhibit double-exponential growth of their vorticity gradients.   We also determine the maximal possible  growth rate and construct solutions that saturate it.   These are the first such results on an unbounded resp.~any 2D domain.
\end{abstract}

\maketitle

\section{Introduction} \lb{S1}

The vorticity form of the incompressible 2D Euler equation modeling ideal fluids is
\beq\lb{1.1}
\omega_t+u\cdot\nabla\omega = 0,  \qquad u:=\nabla^\perp \Delta^{-1} \omega,
\eeq
where $\Delta$ is the Dirichlet Laplacian on the relevant domain, $\nabla^\perp :=(\partial_{x_2},-\partial_{x_1})$, the divergence-free vector field $u$ is the fluid velocity, and  $\omega:=\nabla^\perp\cdot  u$ is its vorticity.  (Our $\nabla^\perp$ is negative of the more standard formula but it will allow us to work with $\omega$ positive in the first quadrant instead of negative.)
On the half-plane $\bbR\times\bbR^+$, the velocity can be found from the vorticity by first extending the latter oddly onto $\bbR\times\bbR^-$ and then convolving it with the Biot-Savart kernel for the Euler equation on $\bbR^2$, namely
\beq\lb{1.2}
u(t,x)=  \frac 1{2\pi} \int_{\bbR^2}  \frac{(x_2-y_2,\,y_1-x_1)}{|x-y|^2} \, \omega(t,y) \, dy.
\eeq

It is a well known result of Yudovich theory that on  smooth 2D domains, initial data from $L^\infty\cap L^1$ give rise to unique global-in-time solutions in this space.  These then satisfy
\[
\|\omega(t,\cdot)\|_{L^p} = \|\omega(0,\cdot)\|_{L^p}.
\]
for all $p\in[1,\infty]$ and $t\ge 0$ because they are transported by measure-preserving flows.  Moreover, Wolibner \cite{Wol} and H\" older \cite{Hol} showed almost a century ago global regularity of solutions in  $W^{1,\infty}\cap L^1$, with  $\|\nabla \omega(t,\cdot)\|_{L^\infty}$ growing no faster than double-exponentially as $t\to\infty$  (even though this rate seems to have first  appeared explicitly only in \cite{Yud}).

Whether this growth rate is in fact attainable has been a major open problem since these works.  In the case of domains without boundaries (i.e., $\bbR^2$ and $\bbT^2$), at least linear growth of gradients was recently proved on $\bbR^2$ for certain perturbations of the Lamb-Chaplygin dipole by Choi and Jeong \cite{ChoJeo}, while the author demonstrated at least exponential growth on $\bbT^2$ \cite{ZlaEulerexp}.  The latter result involves solutions that are $C^{1,\alpha}$ (for any $\alpha\in(0,1)$) and smooth away from a single spatial point. The best result to date for smooth solutions on $\bbT^2$ is an earlier proof of existence of super-linear growth by Denisov \cite{Den}, who also showed that double-exponential growth is possible there on arbitrarily long but finite time intervals \cite{Den2}.

On domains with boundaries (and with the no flow boundary condition),  Yudovich \cite{Yud1, Yud2} and Nadirashvili \cite{Nad} showed prior to the four papers above that unbounded resp.~at least linear growth are possible.
Then, in a striking advance, Kiselev and \v Sver\' ak demonstrated  that infinite time double-exponential growth of the vorticity gradient can indeed occur on a disc \cite{KS}.   Their work was  motivated by numerical simulations of Luo and Hou that suggested existence of  singularity formation for the 3D Euler equation on bounded domains \cite{LuoHou}.

While the result from \cite{KS} was extended to other bounded 2D domains with an axis of symmetry by Xu \cite{Xu}, boundedness of the domain is crucial for the argument to work.  It is based on existence of certain special stable stationary solutions whose perturbations then remain near the stationary solution, and thus the flow velocity  is always a small perturbation  (in some sense) of a specific time-independent background vector field.  The latter exhibits  exponential approach of two particle trajectories to a hyperbolic equilibrium point, and one can then show that some convergent dynamic persists even after perturbation.  The resulting error downgraded the growth to only super-linear in \cite{Den}, while the authors of  \cite{KS} were able to construct perturbations that even enhanced the growth to double-exponential (at the domain boundary, which was a crucial element).

The properties of domains that make existence of stable stationary solutions with the proper dynamic possible in \cite{Den, KS, ZlaEulerexp} do not seem to extend to domains with  infinite areas, even those with a boundary.  In particular, any initial coherence of some non-stationary integrable $\omega$ might be lost after a sufficiently long time because particle trajectories may travel unhindered and possibly escape to infinity.  For this reason, while  the result from \cite{ChoJeo} also holds on the half-plane and thus shows at least linear growth of solution gradients there, no faster growth seems to 
have been known on unbounded domains despite \cite{KS} having been written over a decade ago.  

Our first main result shows that, in fact, double-exponential (i.e.~{\it qualitatively} maximal) gradient growth of smooth solutions does occur in this setting.

\begin{theorem}\lb{T.1.1}
There is $\omega(0,\cdot)\in C_c^{\infty}(\bbR\times\bbR^+)$ with $\|\omega(0,\cdot)\|_{L^\infty}= 1$  such that  the corresponding solution to \eqref{1.1} satisfies
\beq\lb{1.4}
\lim_{t\to\infty} \frac{ \ln\ln \|\nabla \omega(t,\cdot)\|_{L^\infty} } t =\frac2\pi.
\eeq
\end{theorem}

{\it Remarks.}  1.  Simple time scaling shows that with $\|\omega(0,\cdot)\|_{L^\infty}= \gamma$ instead, for some $\gamma>0$, one can obtain \eqref{1.4} with right-hand side $\frac{2\gamma}\pi$.
\smallskip

2.  Note that \eqref{1.4} involves a limit rather than $\liminf$, so the growth rate is exact.
\smallskip

3. Our proof shows that there is $\eps>0$ such that \eqref{1.4} holds whenever  $\omega (0,\cdot)$ is odd in $x_1$ and $\chi_{(\eps,1)\times(0,1)}\le \omega (0,\cdot) \chi_{(0,\infty)^2} \le \chi_{(0,\infty)^2}$.  Scaling shows that here the rectangle $(\eps,1)\times(0,1)$ can be replaced by any multiple of itself.
\medskip

Unlike \cite{Den,KS, ZlaEulerexp}, our proof does not employ any background flow and uses very limited control of most of the mass of the solution.  We only track the motion of a pair of small regions whose areas decay at least exponentially in time, and possibly even double-exponentially (regions of this type were the second main ingredient in \cite{KS}, which also inspired our approach).  Nevertheless, by carefully constructing these evolving regions, we are able to show that they can propel themselves towards each other double-exponentially fast, without contribution from  and  despite possible interference from the rest of the solution.  These regions have different values of the vorticity, they are touching the boundary of the half-plane --- which is crucial just as in \cite{KS} --- and they converge at the origin as $t\to\infty$.  

To obtain the largest possible growth rate in \eqref{1.4}, we construct these regions in such a way that asymptotically they completely fill  a (double-exponentially) shrinking neighborhood of the origin.  We do this because after answering the double-exponential growth question in the affirmative, the next one concerns the maximal rate of such growth.  Since multiplying the vorticity by a constant speeds up the dynamic (and hence increases the rate) by the same factor, this becomes a  question of the maximal possible growth rate for solutions with the constraint $\|\omega(0,\cdot)\|_{L^\infty}= 1$.  This has remained open to date on any domain. 

Our second result implies, together with Theorem \ref{T.1.1}, that the rate $\frac 2\pi$ in \eqref{1.4} is 
 also the maximal possible one on the half-plane.  So gradient growth of the solutions in Theorem~\ref{T.1.1} is also {\it quantitatively} maximal, and the constant $\frac 2\pi$ in \eqref{1.3} below is optimal on $\bbR\times\bbR^+$. 

\begin{theorem}\lb{T.1.2}
For any $\omega(0,\cdot)\in (W^{1,\infty}\cap L^1)(\bbR\times\bbR^+)$ or $\omega(0,\cdot)\in (W^{1,\infty}\cap L^1)(\bbR^2)$ with $\|\omega(0,\cdot)\|_{L^\infty}= 1$, the corresponding solution to \eqref{1.1} satisfies
\beq\lb{1.3}
\limsup_{t\to\infty} \frac{ \ln\ln \|\nabla \omega(t,\cdot)\|_{L^\infty} } t\le \frac 2\pi.
\eeq
\end{theorem}

Proofs of the above results essentially only need to control solutions in small and fast-shrinking neighborhoods of individual points (just the origin in Theorem \ref{T.1.1}).   They can therefore be easily adapted to other smooth domains (with a boundary and an axis of symmetry crossing it in Theorem \ref{T.1.1}) via precise estimates on  the corresponding Dirichlet Green's functions  \cite{GT,Xu} (particularly arbitrarily close to their poles), thus yielding the following.

\begin{theorem}\lb{T.1.3}
Theorem \ref{T.1.2} extends to any uniformly smooth domain in $\bbR^2$, and so does Theorem~\ref{T.1.1} if that domain also has an axis of symmetry that crosses its boundary.
\end{theorem}

{\it Remark.}
In particular, the rate $\frac 2\pi$ is again maximal on any domain of the latter type.
\smallskip

We prove Theorem \ref{T.1.2} in Section \ref{S2} and Theorem \ref{T.1.1} in Section \ref{S3}.  We leave adjustments yielding the proof of Theorem \ref{T.1.3}, following the approach in \cite{Xu},  to the reader.

\medskip

{\bf Acknowledgement.}  
The author was supported in part by NSF grant DMS-2407615.

\section{Proof of Theorem \ref{T.1.2}} \lb{S2}

We start with a crucial estimate for $\omega$ odd in $x_1$, which is a version of results from \cite{KS,ZlaEulerexp}.
For any $x\in[0,\infty)^2$ we let $\til x:=(-x_1,x_2)$, $\bar x:=(x_1,-x_2)$, and 
$Q(r):=(0,\infty)^2\setminus \overline{B_{r}(0)}$.
If $\omega \in (L^\infty\cap L^1)(\bbR\times\bbR^+)$ is odd in $x_1$ (we drop $t$ from the notation in most of our analysis), then the corresponding velocity $u$ is given by
\beq \lb{2.0}
\begin{split}
u_1(x) 
& =  \frac {2x_1}{\pi} \int_{(0,\infty)^2}  \left( \frac{y_1(x_2-y_2)}{|x-y|^2 |x-\til y|^2}
-\frac{y_1(x_2+y_2)}{|x-\bar y|^2 |x+ y|^2}  \right) \omega(y) \, dy, \\
u_2(x) 
& =  - \frac {2x_2}{\pi} \int_{(0,\infty)^2}  \left( \frac{y_2(x_1-y_1)}{|x-y|^2 |x-\bar y|^2}
-\frac{y_2(x_1+y_1)}{|x-\til y|^2 |x+ y|^2}  \right) \omega(y) \, dy.
\end{split}
\eeq
The leading term in both kernels as $|y|\to\infty$ is $-\frac{2y_1y_2}{|y|^4}$, and it turns out that after subtracting it, both integrals in \eqref{2.0} become bounded in all $x$ for which $\frac{x_2}{x_1}$ stays away from 0 and $\pm\infty$.

\begin{lemma}\lb{L.2.1}
There is $C\ge 1$ such that if  $\omega \in (L^\infty\cap L^1)(\bbR\times\bbR^+)$  is odd in $x_1$, then for each $x\in[0,\infty)^2\setminus\{(0,0)\}$ we have
\beq\lb{2.1}
u_j(x)=   (-1)^j \left( \frac 4\pi\int_{Q(|x|)}\frac{y_1y_2}{|y|^4} \omega(y) \, dy 
+  b_j(x)    \right) x_j  \qquad (j=1,2),
\eeq
where $b_1,b_2$ depend on $\omega$ and satisfy
\begin{align*}
|b_1(x)|\le C \|\omega\|_{L^\infty} \left( 1 + \min \left\{ \ln  \frac {x_1+x_2}{x_1}  , \frac{\|\nabla\omega\|_{L^\infty((0,2x_2)^2)}}{\|\omega\|_{L^\infty}} \, x_2 \right\} \right), 
\\ |b_2(x)|\le C\|\omega\|_{L^\infty} \left( 1 + \min \left\{ \ln  \frac {x_1+x_2}{x_2} , \frac{\|\nabla\omega\|_{L^\infty((0,2x_1)^2)}}{\|\omega\|_{L^\infty}}\,  x_1  \right\} \right). 
\end{align*}  
\end{lemma}

{\it Remark.}  
A similar result was proved on the torus in \cite{ZlaEulerexp}.  We will never use the second terms in both $\min$ here, and dropping them yields a result similar to one on the disc from \cite{KS}.  
Our version mainly differs from these two results in the integration domain $Q(|x|)$, a change that --- while technically trivial --- allows for a better comparison of $u$ at different points.

\begin{proof}
First note that the claim of the lemma is equivalent to the one obtained when we replace $Q(|x|)$  by $Q'(2x):=(2x_1,\infty)\times(2x_2,\infty)$.  This is because the difference of these two domains is contained in $D:=[(0,2|x|)\times(\frac{|x|}2,\infty)]\cup [(\frac{|x|}2,\infty)\times(0,2|x|)]$, and
\[
\int_{D}\frac{y_1y_2}{|y|^4}  dy \le 2\int_{(0,2|x|)\times({|x|}/2,\infty)}\frac{y_1}{|y|^3}dy \le 8.
\]

This result with $Q''(2x):=(2x_1,1)\times(2x_2,1)$ in place of $Q'(2x)$ is proved in \cite{ZlaEulerexp} for bounded solutions on $(2\bbT)^2=[-1,1)^2$ (with $-1$ and 1 identified) that are odd in both $x_1$ and $x_2$, and for $x\in[0,\frac 12]^2\setminus\{(0,0)\}$.   This is the same as considering bounded solutions on $\bbR^2$ that are odd and 2-periodic in both variables, with all integrals being principal value ones.  We  consider here $L^\infty\cap L^1$ solutions on $\bbR^2$ that are odd in both variables (after an odd extension of $\omega$ to $\bbR^2$), so we can drop from the proof in \cite{ZlaEulerexp} all estimates involving integrals over $\bbR^2\setminus(-1,1)^2$, and  perform the estimates involving integrals over $(-1,1)^2$ instead with integrals over $\bbR^2$.  These all become appropriate integrals over $(0,1)^2$ resp.~$(0,\infty)^2$ due to the symmetries of $\omega$, as in \eqref{2.0}, and  the argument from \cite{ZlaEulerexp} now applies after simply changing all upper bounds 1 in the various integrals to $\infty$ (and $Q''(2x)$ to $Q'(2x)$).  We refer the reader to 
the short proof of \cite[Lemma 2.1]{ZlaEulerexp} 
for details.
\end{proof}

Consider now points $(\pm\eps,0)\in\bbR^2$ with $\eps\in(0,1]$ and any $\omega \in (L^\infty\cap L^1)(\bbR^2)$ with $\|\omega\|_{L^\infty}=1$ and $\|\omega\|_{L^1}\le I$ for some fixed $I\ge 1$.  From \eqref{1.2} we see that the instantaneous approach velocity of the points $(\pm\eps,0)$  is
\[
u_1(-\eps,0)-u_1(\eps,0) = \frac {2\eps}{\pi} \int_{\bbR^2}  \frac{y_1y_2}{|(\eps,0)-y|^2 \,|(-\eps,0)-y|^2} \, \omega(y)dy, 
\]
which is clearly maximized when $\omega(y)=\sgn (y_1y_2)$ on the super-level set 
\[
L_\eps(a_{I,\eps}):=\left\{ y\in\bbR^2 \,\bigg|\, \frac{|y_1y_2|}{|(\eps,0)-y|^2 \,|(-\eps,0)-y|^2} \ge a_{I,\eps} \right\} 
\]
with $a_{I,\eps}>0$ such that $|L_\eps(a_{I,\eps})|=I$, and $\omega=0$ on $\bbR^2\setminus L_\eps(a_{I,\eps})$.  Since $L_\eps(a_{I,\eps})$ is symmetric in both coordinates, we get
\beq\lb{2.20}
u_1(-\eps,0)-u_1(\eps,0) \le  \frac {8\eps}{\pi} \int_{(0,\infty)^2\cap L_\eps(a_{I,\eps})}  \frac{y_1y_2}{|(\eps,0)-y|^2 \,|(-\eps,0)-y|^2} \, dy, 
\eeq
with equality holding in the case described above.  But this maximal case is also realized by an $\omega\in (L^\infty\cap L^1)(\bbR\times\bbR^+)$ that is odd in $x_1$ (namely with $\omega\chi_{(0,\infty)^2}= \chi_{(0,\infty)^2\cap L_\eps(a_{I,\eps})}$), so Lemma~\ref{L.2.1} and oddness of $u_1$ in $x_1$ show that then
\[
\left| u_1(-\eps,0)-u_1(\eps,0)  - \frac {8\eps}\pi\int_{Q(\eps)\cap L_\eps(a_{I,\eps})}\frac{y_1y_2}{|y|^4} \,  dy 
  \right| \le 2 C\eps .
\]

For all  $\eps\in(0,1]$ we  have $L_\eps(a_{I,\eps})\subseteq B_{10\sqrt I}(0)$ because $I\ge 1$, 
so we always obtain
\begin{align*}
u_1(-\eps,0)-u_1(\eps,0)  & \le  \frac {8\eps}\pi\int_{Q(\eps)\cap B_{10\sqrt I}(0)}\frac{y_1y_2}{|y|^4} \,  dy 
+ 2C\eps 
\\ & = \frac {8\eps}\pi   \ln\frac{10\sqrt I}\eps \int_0^{\pi/2}  \cos\theta \sin\theta \, d\theta +2C\eps
\\ & =  \frac {-4\ln\eps  + 4\ln(10\sqrt I)+2C\pi}\pi \,\eps.
\end{align*}

However, if also $ \|\nabla \omega\|_{L^\infty} \le \frac 1\eps$,
then we can use $|\omega(y_1,y_2)-\omega(-y_1,y_2)|\le 2\|\nabla \omega\|_{L^\infty}|y_1|$ when $|y_1|\le \|\nabla \omega\|_{L^\infty}^{-1}$ (and cancellation of the kernel at these two points) in \eqref{2.20} to refine  the above estimates and obtain with $m:=\min\{\|\nabla \omega\|_{L^\infty}^{-1}, 1\}$,
\begin{align*}
u_1(-\eps,0)-u_1(\eps,0)  & \le  \frac {8\eps}\pi\int_{Q(\eps)\cap B_{10\sqrt I}(0)}\frac{y_1y_2}{|y|^4} \, \min\{\|\nabla \omega\|_{L^\infty}|y_1|,1\} \, dy 
+ 2C\eps 
\\ & \le   \frac {8\eps}\pi\int_{Q(m)\cap B_{10\sqrt I}(0)}\frac{y_1y_2}{|y|^4} \, dy 
+  \frac {8\eps \|\nabla \omega\|_{L^\infty}}\pi\int_{Q(\eps)\setminus Q(m)}\frac{dy}{|y|} 
+ 2C\eps 
\\ & \le   \frac {-4\ln m +4\ln(10\sqrt I)+(4+2C)\pi}\pi \,\eps.
\end{align*}

Rotational symmetry of the equation on $\bbR^2$ now shows that any two points with distance $\delta:=2\eps\in(0,2]$ have instantaneous approach velocity no greater than  
\beq\lb{2.30}
\frac {2 \min \left\{-\ln \frac\delta 2, - \ln m \right\}+2C_I}\pi \,\delta = \frac 2\pi \left( \ln \min\left\{2\delta^{-1}, \max\left\{ \|\nabla \omega\|_{L^\infty},1 \right\} \right\}+C_I \right)  \delta,
\eeq
 with $C_I:= \ln(10\sqrt I)+(1+\frac C2)\pi$.   This means that 
 \[
 \frac d{dt} \|\nabla \omega(t,\cdot)\|_{L^\infty} \le \frac 2\pi  \big(\ln \max\{\|\nabla \omega(t,\cdot)\|_{L^\infty},1\} +C_I\big) \, \|\nabla \omega(t,\cdot)\|_{L^\infty}
 \]
 holds for solutions on $\bbR^2$ with  $\|\omega(0,\cdot)\|_{L^\infty}=1$ and $\|\omega(0,\cdot)\|_{L^1}\le I$ for all $t\ge 0$.

Since any solution on $\bbR\times\bbR^+$ also becomes a solution on $\bbR^2$ after odd reflection in $x_2$, the same bound holds for the former (when $\|\omega(0,\cdot)\|_{L^\infty}=1$ and $\|\omega(0,\cdot)\|_{L^1}\le \frac I2$).  Hence 
\[
 \|\nabla \omega(t,\cdot)\|_{L^\infty} \le \exp\left(   e^{2t/\pi+C_I(1-e^{-2t/\pi}) } \max  \left\{\ln \|\nabla \omega(0,\cdot)\|_{L^\infty}, 1 \right\} \right)
\]
holds on both domains for all $t\ge 0$, and \eqref{1.3}  follows.

\section{Proof of Theorem \ref{T.1.1}} \lb{S3}

Let
\[
g(s):=
se^{|\ln s|^{1/2}} 
\]
for $s\ge 0$, and  for any $\eps\in[0,e^{-1})$ let
\begin{align*}
\Omega_\eps := \big\{ x\in(0,\infty)^2\,\big|\, x_1\in(\eps,e^{-1}) \,\,\&\,\, x_2<g(x_1) \big \}.
\end{align*}
Then for all $s\in(0,e^{-4}]$ and $D_s:=\Omega_0\cap Q\big(\sqrt{s^2+g(s)^2}\,\big)\cap B_{e^{-1}}(0)$ we have
\beq \lb{6.1}
\frac 4\pi  \int_{D_s} \frac{y_1y_2}{|y|^4} dy =  h(s) |\ln s| ,
\eeq
for some $h:(0,e^{-4}]\to (0,\infty)$ such that 
\beq\lb{6.2}
\lim_{s\to 0} h(s)=\frac 2\pi.
\eeq
Indeed, these claims follow from
\[
\lim_{s\to 0} \frac{|\ln g(s)|}{|\ln s|}= \lim_{s\to 0} (1-|\ln s|^{-1/2}) =1,
\]
$\lim_{s\to 0} \frac{g(s)}s=\infty$, and
\[
\frac 4\pi  \int_0^{\pi/2}  \cos\theta \sin\theta \, d\theta = \frac 2\pi.
\]

Let $C$ be from Lemma \ref{L.2.1} and let $s_0\in(0,e^{-4}]$ be such that with the function
\beq\lb{6.11}
f(s):=2+ |\ln s|^{1/2} \ge 1+\ln\frac{s+g(s)}{s}
\eeq
(recall that $\ln (1+\beta ) \le 1+ \ln \beta$ for $\beta\ge 1$) we have
\[
h(s)|\ln s|\ge Cf(s) 
\]
 for all $s\in(0,s_0]$.  Also let
\[
C':=C\, \frac{g(s_0)}{s_0} \left( 1+\ln \frac{s_0+g(s_0)}{s_0} \right) 
= \max_{s\in[s_0,e^{-1}]} C\, \frac{g(s)}{s} \left( 1+\ln \frac{s+g(s)}{s} \right) .
\]

Next note that the unit outer normal vector to  $\Omega_\eps$ at $(s,g(s))$ is
$v_s(-1,\frac 1{g'(s)})$, where 
\[
v_s:=\left( 1+g'(s)^{-2} \right)^{-1/2} \, ,
\]
and for $s\in(0,e^{-1}]$ we have
\[
g'(s)= e^{|\ln s|^{1/2}} \left(1-\frac 1{2|\ln s|^{1/2}} \right) \in\left(1, \frac{g(s)}s \right).
\]
We then let
\beq\lb{6.10}
\rho_0:= \frac{s_0}{2 |\ln s_0|^{1/2}} \le  s_0 \left( \left(1-\frac 1{2|\ln s_0|^{1/2}} \right)^{-1} -1\right) = \min_{s\in[s_0,e^{-1}]} \left( \frac{g(s)}{g'(s)}-s \right) . 
\eeq

If now $\eps\in(0,s_0]$, and  $\omega \in (L^\infty\cap L^1)(\bbR\times\bbR^+)$ is odd in $x_1$ and $\chi_{\Omega_\eps}\le \omega \chi_{(0,\infty)^2} \le \chi_{(0,\infty)^2}$, Lemma \ref{L.2.1}, our choice of $s_0$, and \eqref{6.1} show that the corresponding velocity $u$ satisfies 
\begin{align*}
(-1)^j u_j(s,g(s))\ge 0   & \qquad\text{for $s\in[\eps,s_0]$ and $j=1,2$}, 
\\ (-1)^j u_j(s,g(s)) \ge -C's  & \qquad\text{for $s\in[s_0,e^{-1}]$ and $j=1,2$}. 
\end{align*}  
Using also \eqref{6.2} and \eqref{6.11} gives us
\begin{align*}
 \eps^{-1} \sup_{s\in[0,g(\eps)]} u_1(\eps,s)   \le  -h(\eps)|\ln\eps|+Cf(\eps)+3C  + e \inf_{s\in[0,1]} u_1(e^{-1},s),
 \end{align*}  
and simple scaling then yields the following result.

\begin{lemma}\lb{L.6.1}
Let $\eps\in(0,s_0]$ and $\alpha>0$.  If  $\omega \in (L^\infty\cap L^1)(\bbR\times\bbR^+)$  is odd in $x_1$ and $\chi_{\alpha \Omega_\eps}\le \omega \chi_{(0,\infty)^2} \le \chi_{(0,\infty)^2}$,  then 
\begin{align}
(-1)^j u_j(s\alpha,g(s)\alpha)\ge 0   & \qquad\text{for $s\in[\eps, s_0]$ and $j=1,2$}, \lb{6.3}
\\ (-1)^j u_j(s\alpha,g(s)\alpha) \ge -C's\alpha  & \qquad\text{for $s\in[s_0,e^{-1}]$ and $j=1,2$}, \lb{6.4}
\end{align}  
as well as
\begin{align}
 (\eps\alpha)^{-1} \sup_{s\in[0,g(\eps)\alpha]} u_1(\eps\alpha,s)   
 \le  -h(\eps)|\ln\eps|+Cf(\eps)+3C  + (e^{-1}\alpha)^{-1} \inf_{s\in[0,\alpha]} u_1(e^{-1}\alpha,s) .   \lb{6.5}
\end{align}  
\end{lemma}

Consider  now some decreasing $\eps:[0,\infty)\to (0,s_0]$ (to be determined),  and let $\omega$ be an odd-in-$x_1$ solution to \eqref{1.1} with  $\omega(0,\cdot) \in C^\infty_c(\bbR\times\bbR^+)$ such that $\chi_{\Omega_{\eps(0)}}\le \omega (0,\cdot)\chi_{(0,\infty)^2} \le \chi_{(0,\infty)^2}$.  Then of course $\|\omega(t,\cdot)\|_{L^\infty}=1$ for all $t\ge 0$.  Also let $\alpha:[0,\infty)\to(0,1]$ satisfy $\alpha(0)=1$ and 
\beq \lb{6.6}
\alpha'(t)=  -3(C's_0\rho_0^{-1}+C)\alpha(t) +  e \inf_{s\in[0,\alpha(t)]} u_1(t,e^{-1}\alpha(t),s)  
\eeq
for $t\ge 0$, with $u(t,\cdot)$ the velocity for $\omega(t,\cdot)$.  Note that Lemma \ref{L.2.1} shows that 
\beq \lb{6.7}
\alpha'(t)\le -3C' s_0\rho_0^{-1}\alpha(t) \qquad (<0),
\eeq
 and  $\alpha(t)$  remains positive because from \eqref{2.30} we see that
\begin{align*}
 |u_1(t,e^{-1}\alpha(t),s)  |  \le   |u_1(t,e^{-1}\alpha(t),s) -  u_1(t,0,s)| \le
\frac {2 |\ln (e^{-1}\alpha(t))| +C_{\frac 12\|\omega(0,\cdot)\|_{L^1}}}{\pi e} \, \alpha(t).
 \end{align*}

We now want to show that the boundary of the level set $\{ \omega(t,\cdot)=1 \}$ can never enter  $\alpha(t)\Omega_{\eps(t)}$ for appropriate functions $\eps$.
At any given time $t\ge 0$, the outer normal velocity of the boundary of the evolving domain $\alpha(t)\Omega_{\eps(t)}$ is $e^{-1}\alpha'(t)$ on  $\{(e^{-1}\alpha(t),s)\,|\, s\in[0, \alpha(t)]\}$, where the unit outer normal vector is $(1,0)$ and 
\[
u(t,e^{-1}\alpha(t),s)\cdot (1,0) = u_1(t,e^{-1}\alpha(t),s) > e^{-1} {\alpha'(t)}
\]
by \eqref{6.6}.  From \eqref{6.10} and \eqref{6.7} we see that that velocity is
\[
\alpha'(t)(s,g(s))\cdot v_s(-1,{g'(s)}^{-1}) 
\begin{cases}
 <0 & s\le s_0,
\\ \le \rho_0 v_s \alpha'(t) \le - 3C's_0 v_s\alpha(t)  & s\ge s_0
\end{cases}
\]
on $\{(s\alpha(t),g(s)\alpha(t))\,|\, s\in[\eps(t),e^{-1}]\}$, where the unit outer normal vector is $v_s(-1,\frac 1{g'(s)})$ and 
\[
u(t,s\alpha(t),g(s)\alpha(t))\cdot v_s(-1,{g'(s)}^{-1}) \ge  0
\]
by \eqref{6.3} when $s\in [\eps(t),s_0]$, while
\[
u(t,s\alpha(t),g(s)\alpha(t))\cdot v_s(-1,{g'(s)}^{-1}) \ge  -2C'v_s \,s\alpha(t) > - 3C's_0 v_s\alpha(t) 
\]
by \eqref{6.4} when $s\in [s_0, e^{-1}]$.
(Note that in both corners we consider vectors normal to both segments of $\alpha(t)\Omega_{\eps(t)}$ meeting at that corner.)   Finally, the velocity is $-(\eps(t)\alpha(t))'$ on $\{(\eps(t)\alpha(t),s)\,|\, s\in[0,g(\eps(t))\alpha(t)]\}$, where the unit outer normal vector is $(-1,0)$ and 
\[
u(t,\eps(t)\alpha(t),s)\cdot (-1,0) 
\ge \big[ h(\eps(t))|\ln\eps(t)|-Cf(\eps(t))-3C \big]  \eps(t)\alpha(t) - e\, {\eps(t)} \inf_{s\in[0,\alpha(t)]} u_1(t,e^{-1}\alpha(t),s)
\]
by \eqref{6.5}.  Recalling \eqref{6.6} and the definition of $f$, the last expression will be strictly greater than $ -(\eps(t)\alpha(t))'$ if 
\beq\lb{6.9}
\eps'(t) > -\big[  h(\eps(t))|\ln\eps(t)| -C|\ln\eps(t)|^{1/2}-8C-3C' s_0\rho_0^{-1} \big] \eps(t)
\eeq
holds for all $t\ge 0$.  If this is satisfied,  the strict inequalities above  mean that the boundary of  $\{ \omega(t,\cdot)=1 \}$ can indeed never enter $\alpha(t)\Omega_{\eps(t)}$.  

So then $\omega(t,\eps(t)\alpha(t),s)=1$ for all $t\ge 0$ and $s\in(0,g(\eps(t))\alpha(t)]$.  Since  $\omega(t,0,s)=0$ for these $s$ by symmetry, it follows that
\beq\lb{6.12}
\|\nabla \omega(t,\cdot)\|_{L^\infty} \ge (\eps(t)\alpha(t))^{-1} \ge \eps(t)^{-1} 
\eeq
 for all $t\ge 0$.  But \eqref{6.2} shows that when $\eps(0)>0$ is small enough, \eqref{6.9} has a decreasing solution satisfying 
 \[
 \lim_{t\to\infty} \frac {\ln(-\ln \eps(t))}t = \frac 2\pi.
 \]
So  \eqref{6.12} yields
 \[
 \liminf_{t\to\infty} \frac{ \ln\ln \|\nabla \omega(t,\cdot)\|_{L^\infty} } t\ge \frac 2\pi,
 \]
 which together with \eqref{1.4} finishes the proof. 

\end{document}